\documentstyle[12pt]{article}
\newcommand{\const}{\mathop{\rm const}\limits}

\newcommand{\supp}{\mathop{\rm supp}\limits}

\begin{document}

\begin{center}

{\bf Trigonometric Approximation in Sobolev-Grand Lebesgue Spaces. }\\

\vspace{3mm}

{\sc Ostrovsky E., Sirota L.}\\

\normalsize

\vspace{3mm}

{\it Department of Mathematics and Statistics, Bar-Ilan University, 59200, \\
Ramat Gan, Israel.}\\
e-mail:\ eugostrovsky@list.ru \\

\vspace{3mm}

{\it Department of Mathematics and Statistics, Bar-Ilan University, 59200, \\
Ramat Gan, Israel.}\\
e-mail:\ sirota3@bezeqint.net \\

\vspace{4mm}

{\bf Abstract.}\\

\end{center}

\vspace{4mm}

  \ We study in this short preprint the theory of trigonometric approximation
in the so-called Banach functional rearrangement invariant Sobolev-Grand Lebesgue Spaces. \par

\vspace{4mm}

 {\it Key words and phrases:}  Rearrangement invariant Banach functional spaces, trigonometric polynomials
 and approximation, Sobolev and Grand Lebesgue Spaces, module of continuity, fundamental function, indicator;
 direct and inverse theorem and estimates;
  Fejer, Jackson, Fourier, Dirichlet and Vallee Poussin  kernels, Kolmogorov's width.\\

{\it Mathematics Subject Classification (2000):} primary 60G17; \ secondary
 60E07; 60G70.\\

\vspace{3mm}

\section{ Introduction. Notations. Statement of problem.}

\vspace{3mm}

 \hspace{4mm} Let $  \ X = [0, 2 \pi] $ with {\it normed} Lebesgue measure $ \ d \mu(x) = \mu(dx) = dx/(2 \pi) \  $ be  a
classical probability space and let $  B  $
be  Banach rearrangement invariant real valued functional space builded over $  X, $ equipped with norm $ ||f||B, $
consisting on the $  \ 2 \pi \ - \ $ periodical functions. \par

 \ The space  $ B $ is called a {\it homogeneous Banach space,} (abbreviated HBS), on the set $ \ X, \  $ see
[18], [10], [29] etc., iff: \\

 \ {\bf (a).} It is linear subspace of $ L_1 = L_1(X, \mu),  $ such that
$$
{\bf (a).} \hspace{6mm} \exists C \in (0, \infty) \ \Rightarrow \ ||f||_1 \le C \ ||f||B, \  \eqno(1.1)
$$

 \ {\bf (b).} The translation $ \ U_t[f], \  U_t[f](x) := f(x - t) $ is continuous isometry of $ \ B \ $ onto itself, namely:

$$
{\bf (b).} \hspace{6mm} \forall f \in B \ \Rightarrow ||U_t f||B = ||f||B \eqno(1.2)
$$
and \ {\bf (c).}

$$
{\bf (c).} \hspace{6mm}  \lim_{t \to 0} || U_t[f] - f ||B = 0. \eqno(1.3)
$$

 \ For instance, the classical Lebesgue - Riesz spaces $ \ L_p = L(p) = L_p(X,\mu), \ 1 \le p < \infty \ $ with ordinary norms

$$
|f|_p = |f|L_p = |f|L(p) = |f|L_p(X,\mu) = \left[ \int_X |f(x)| \ \mu(dx)  \right]^{1/p}
$$
are  HBS. \par

 \ Define for each  such a function $  f = f(x), \ f: X \to R $ from this space $ \ B: \ f \in B \ $ its $ \ B \ -  $ module of continuity

$$
\omega_B[f](\delta) \stackrel{def}{=} \sup_{h: |h| \le \delta} ||U_h[f] - f||B, \ \delta \in X. \eqno(1.4)
$$

 \ We admit in (1.4) that $ \ f(x) = 0 $ if $ x \notin X, \ $  and will write for brevity

$$
\omega[f](\delta)_p = \omega_{L(p)}[f](\delta)=
 \sup_{h: |h| \le \delta} |U_h[f] - f|_p, \ \delta \in X. \eqno(1.4a)
$$

 \ Evidently,

$$
B \in HBS, \ f \in B \ \Rightarrow \ \lim_{\delta \to 0+}\omega_B[f](\delta)  = 0.
$$

\vspace{3mm}

 \ Note that in the approximation theory may be successfully  used another modules of continuity, see [10].
 The case of general rearrangement spaces $ \ B \ $ is studied in [4], [21], [31], [32]. \par

\vspace{4mm}

 \ {\bf The purpose of this note is to extend the characterization of best trigonometric
approximation by the module $ \omega_B[f](\delta)  $  to the suitable subspace of the so-called Banach space of periodic functions,
namely, Grand Lebesgue Spaces (GLS), as well as to the so-called exponential Orlicz spaces.\par

 \ The case when the space $  B  $ coincides with some Sobolev space $ \ W_p^r, \ p \ge 1, \ r = 0,1,2,\ldots  $  will be also
 considered further. } \par

\vspace{4mm}

 \ Hereafter $ C, C_j $ will denote any non-essential finite positive constants. \par
 \ Further, let $ (Y, ||\cdot||Y) $ be any rearrangement invariant (r.i.) space on the set $ X; $ denote by $ \phi(Y, \delta) $
  its fundamental function

  $$
  \phi(Y,\delta) = \sup_{A, \mu(A) \le \delta} ||I(A)||Y, \eqno(1.5)
  $$
  where as usually  $ \ I(A) \ $ denotes the ordinary indicator function of the measurable set $ \ A: $

$$
\ I(A)= I(A,x) =   1, x \in A; \  I(A) = I(A,x) = 0, \ x \notin  A.
$$

\vspace{3mm}

 \ Denote by $ T(n), \ n = 1,2,\ldots $ the set (subspace) of all the trigonometric polynomials on the variable $ x; \ x \in X \ $ of degree
less or equal $  n $ with coefficients from the space $  B  $
and  define correspondingly for the space $  B, $ in particular,  for the space $  \ G(\psi), $ and for
each function $ \ f $ from this space the minimal error of its trigonometrical approximation

$$
E_n[f]B \stackrel{def}{=} \inf_{g \in T(n)} ||f - g||B. \eqno(1.6)
$$

\ We have to take as above in (1.6) for brevity

$$
E_n[f]_p \stackrel{def}{=} E_n[f]L_p = \inf_{g \in T(n)} |f - g|_p. \eqno(1.7)
$$

\vspace{3mm}

\ {\bf Definition 1.1.} The function $  f  $ from the Banach space   $  B  $ is said to be  {\it trigonometric approximated}
in this space, write: $  \ f \in TA(B), \ $ if

$$
\lim_{n \to \infty} E_n[f]B = 0. \eqno(1.8)
$$

\vspace{3mm}

\section{ The case of Grand Lebesgue Spaces.}

\vspace{3mm}

 \hspace{6mm} We recall first of all some needed facts about  Grand Lebesgue Spaces (GLS). \par

\ Recently, see [5], [11], [12], [15]-[17], [20], [22]-[27]
etc.  appear the so-called Grand Lebesque Spaces $ GLS = G(\psi) =
G(\psi; b), \ b = \const \in (1, \infty] $ spaces consisting on all the measurable functions $ f: X \to R $ with finite norms

$$
     ||f||G(\psi) = G(\psi,b) \stackrel{def}{=} \sup_{p \in [1,b)} \left[ \frac{|f|_p}{\psi(p)} \right]. \eqno(2.1)
$$

  \ Here $ \psi(\cdot) $ is some continuous positive on the semi - open interval $ [1,b) $ function such that

$$
     \inf_{p \in [1,b)} \psi(p) > 0.
$$

 \ It is evident that $ G(\psi; b) $ is Banach functional rearrangement invariant (r.i.) space and
 $ \supp(G(\psi_b)) := \supp \psi_b = [1,b). $\par

\ Let the {\it family} of measurable functions
$ h_{\alpha} = h_{\alpha}(x), \ x \in X, \ \alpha \in A,  $ where $  \ A \ $ be arbitrary set, be such that

$$
\exists b \in (1, \infty], \ \forall p \in [1, b) \ \Rightarrow
\psi^A(p) := \sup_{\alpha \in A} | \ h_{\alpha} \ |_p < \infty.
$$

 \ Such a function $  \psi^A(p) $ is named as a {\it  natural function} for the family $  A. $ Obviously,

$$
\sup_{\alpha \in A} || \ h_{\alpha} \ ||G\psi^A = 1.
$$

 \ Note that the case when  $ \ \sup_{p \in [1,b)} \psi(p) < \infty  $ is trivial for us; if for instance $  \ b < \infty $ and
 $ \psi(b-0) < \infty, $ then the space $ G(\psi,b) $ coincides with ordinary Lebesgue-Riesz space $  L_b(X). $ Therefore,
 we can and will suppose in the sequel without loss of generality

 $$
  \sup_{p \in [1,b)} \psi(p) = \lim_{p \to b-0} \psi(p) = \infty. \eqno(2.2)
 $$

  \ These spaces are used, for example, in the theory of probability, theory of PDE, functional analysis
 theory of Fourier series, theory of martingales etc.\par

 \ The problem of trigonometric approximation in the Grand Lebesgue Spaces is considered in more complicated terms in the article [7].\par

 \ The spaces $ G(\psi; b) $ are non-separable, but they satisfy the Fatou property. As long as its
Boyd's indices  $ \gamma_-, \gamma_+ $ are correspondingly

$$
\gamma_- = 1/b, \ \gamma_+ = 1,
$$
 we conclude  that the spaces $ G(\psi; b) $ are interpolation spaces not only between the spaces
$ [L_1, \ L_{\infty}] $ but between also the spaces $ [L_1, \ L_s] $ for all values $ s $ for which
$ \ s > b. $ \par

\vspace{3mm}

 \ Note that an another approach to the problem of combine the $ L_p $ spaces may be found in the articles [2], [3],
[8], [30]. \par

 \vspace{3mm}

 \ Let us introduce an important subspace of the whole space $ G(\psi; b).  $ \par

\vspace{3mm}

 \ {\bf Definition 2.1.} The (closed) subspace $ G^o(\psi) = G^o(\psi; b)  $ of the whole GLS $ \ G\psi: \ G^o(\psi) \subset G(\psi) \  $
 consists by definition on all the functions $ \ \{f\} \ $ from the whole space $ G(\psi; b)  $ for which

$$
\lim_{p \to b - 0} \left\{ \frac{|f|_p}{\psi(p)} \right\} = 0.\eqno(2.3)
$$

 \ Of course, the functions belonging to the space $ G^o(\psi) = G^o(\psi; b)  $ have at the same norm (2.1) as in the space
 $ G(\psi) = G(\psi; b). \  $  \par

 \ It is known, see [26], [43], that the spaces $ G^o(\psi) = G^o(\psi; b)  $ have absolute continuous norm and coincides with closure of
the set of all bounded measurable functions. Alike in the theory of Orlicz spaces, they are reflexive and separable.\par

\vspace{3mm}

 \ {\bf Example 2.1.} Let the numerical valued {\it random variable (r.v.)} (measurable function) $ \xi = \xi(x), \ x \in X  $
be defined on our probability space and has a standard Gaussian (normal) distribution. Then it belongs to the GLS
$  \ G\psi_{1/2},   $ where for each $ \ m = \const \in (0, \infty) $ we define

$$
\psi_{m}(p) := p^{1/m}, \ p \ge 1.
$$

 \ Indeed, it is easily to calculate

$$
|\xi|_p \asymp p^{1/2}, \ p \ge 1.
$$
 \ Therefore, $ \xi \in G\psi_{1/2}, $ but $ \xi \notin G^o\psi_{1/2}. $\par
 \ On the other hand,

$$
 \forall m \in (0,2) \ \Rightarrow \xi \in G^o\psi_{1/m}.
$$

\vspace{4mm}

 \ Let us now consider a generalizations of the classical results of trigonometric approximation into the Grand Lebesgue Spaces.\par

\vspace{3mm}

  \ {\bf Theorem 2.1.} Suppose the function $  f  $ belongs to the space $  G\psi. $ Statement: this function is
 trigonometric approximated in this space,  $  \ f \in TA(G\psi), $ if  and only if it belongs to the subspace
 $ \ G^o(\psi) = G^o(\psi; b).  $  \par
  \ Furthermore, in this case

$$
E_n[f]G\psi \le C_1(G\psi) \ \omega_{G\psi}[f](2 \pi/n), \ n = 1,2,\ldots; \eqno(2.4)
$$
and conversely

 $$
 \omega_{G\psi}[f](2 \pi/n) \le C_2 \ n^{-1} \sum_{k=1}^n E_k[f]G\psi. \eqno(2.5)
 $$

\vspace{3mm}

 \ {\bf Proof.} \par

 \vspace{3mm}

 {\bf 1.} Assume at first  $ \ f \in G^o(\psi), \ $ then the function $ \ f(\cdot) \ $  has an Absolute Continuous Norm (ACN)
 in the space $ \ G\psi $ or equally in the space $ \ G^o\psi, \ $ see [22], [24], [25], [26].
 Following, $ \ G^o(\psi) \ $ is homogeneous Banach space (HBS), see [4], chapter 1, section 3; and hence

$$
\lim_{n \to \infty} \omega_{G\psi}[f](2 \pi/n) = 0. \eqno(2.6)
$$

 \ It remains to apply the classical results about trigonometric approximation in these spaces, see [10], [18], [29] etc.
to deduce the equalities (2.4) and (2.5). \par

\vspace{3mm}

 \ Note in addition to this pilcrow that in (2.4) the apparatus for correspondent approximation can be the convolution
with the classical trigonometric kernels: Fejer $  F_n,$  Jackson $  J_n, $ Fourier $ S_n, $ Riesz $  R_n, $ Dirichlet $ D_n, $
Vallee Poussin $  V_n $ etc. \par

\vspace{3mm}

 \ {\bf 2.}  \ Two examples of the direct estimates. Assume $  f \in G^o(\psi). $ \par

\vspace{3mm}

  \ Let us write the famous Jackson's proposition for the $ L_p, \ p \ge 1 $  spaces:

$$
| \ f - J_n*f \ |_p \le C_1 \ \omega[f](2 \pi /n)_p. \eqno(2.7)
$$

 \ It is important to note that the constant $  C_1  $ one can choose not depending on the  variable $ p. $
Furthermore, the estimate (2.7) there holds still for the value $  \ p = \infty, \ $  where $  C_1 $ may be taken such that
$  C_1 = 3. $\par

 \ As long as $  f \in G\psi,  $

$$
\omega[f](2 \pi /n)_p \le  \omega_{G\psi}[f](2 \pi /n) \ \psi(p), \ p \in (1,b),
$$
and we get after substituting into (2.7)

$$
| \ f - J_n*f \ |_p \le C_1 \ \omega_{G\psi}[f](2 \pi /n) \ \psi(p),
$$
or equally on the basis of the direct definition of $  \ G\psi \ $ norm

$$
|| \ f - J_n*f \ ||G\psi \le C_1 \ \omega_{G\psi}[f](2 \pi /n),
$$
which coincides  with (2.4); the relation (2.5) may be grounded analogously. \par

 \ A second example. Let for definiteness the number $  \ n \ $ be even number, and let  $  f \in G^o\psi =  G^o\psi_b,
 b = \const \in (1,\infty].  $ One can use the-known Vallee-Poussin inequality

$$
| \ f - V_n*f \ |_p \le C_2 \ E_{n/2}[f]_p, 1 \le p < b.
$$

 \ Since

$$
E_{n/2}[f]_p = \sup_{g \in T(n/2)}| \ f - g \ |_p \le \sup_{g \in T(n/2)} || \ f - g \ ||G\psi \cdot \psi(p),
$$
we deduce for any function $  \  f  \ $ from the space $ \ G^o\psi $

$$
 || \ f  - V_{n}*f \ ||G\psi \le C_2 \ \sup_{g \in T(n/2)} || \ f - g \ ||G\psi \  = E_{n/2}[f]G\psi.\eqno(2.8)
$$

 \  Note that the last  expression tends to zero as $ \ n \to \infty \ $  as long as $ f \in \ G^o\psi. $

\vspace{3mm}

{\bf 3.} We continue. Suppose for certain function $ \ f \ $ from the space $  \ G\psi \ $

$$
\lim_{n  \to \infty} E_n[f]G\psi = 0.
$$

\vspace{3mm}

 \ Since the trigonometric system consists only on bounded functions,
  the last equality implies that the function $ \ f  \ $ belongs to the closure in $ \ G\psi \ $ norm of the set of all
 bounded measurable functions $ G^{(b)}G\psi. $ But the last space coincides with the space $ \ G^o\psi. \ $ \par
 \ The rest follows from the theory of r.i. spaces, see the classical book of C.Bennet, R.Sharpley \ [4], chapter 1. \par

\vspace{3mm}

 \ {\bf Definition 2.1.}  Let $  G\psi $ and $  G\nu $ be two Grand Lebesgue Spaces with at the same support of the
generating functions $  \psi $  and $  \nu.$  We will  write  $  \nu  << \psi, $ or equally
$  G\nu  << G\psi, $ iff

$$
\lim_{p \to b - 0} \frac{\psi(p)}{\nu(p)} = 0.\eqno(2.9)
$$

\vspace{3mm}

 {\bf Corollary 2.1.}  Let $  0 \ne f \in G\psi  \setminus G^o\psi  $  and let $ \nu = \nu(p) $ be arbitrary another
$  \ \Psi \ $ function with at the same support as the source function $ \ \psi = \psi(p), \ $  and such that
$  G\nu  << G\psi. $ \ Then the function $ f(\cdot) $ is not trigonometric approximated in the space
$ \ G\psi: \ f \notin TA(G\psi) $  but it is trigonometrical approximated in the space $ \ G\nu: \ f \in TA(G\nu). $ \par

 \vspace{4mm}

\section{ Main result: Approximation in Sobolev-Grand Lebesgue Spaces. }

\vspace{3mm}

 \hspace{4mm} Denote as ordinary by $ W^r_p, \ 1 \le p < \infty, \ r = 1,2,\ldots  $ the classical Sobolev's space on the unit
circle $ X = [0, 2 \pi] $ consisting on the $ \ 2 \pi \ $ periodical  functions $ \ \{f\}, \ $ in particular $ f(0) = f(2 \pi). $ \par

 \ The space $ W^0_p $ coincides with the classical Lebesgue - Riesz space $  L_p(X). $ \par
 \ The norm of a function $ \ f \ $ in this space  $ W^r_p $ may be defined for instance as follows:

$$
||f||W^r_p \stackrel{def}{=} \left[ |f|_p^p +   |f^{(r)}|_p^p \right]^{1/p}. \eqno(3.1)
$$

 \ Let also $ \ \psi = \psi(p), \ 1 \le p < b, \ b = \const \in (1, \infty] $ be the ordinary $ \ \Psi \ -  $ function. \par

\vspace{4mm}

 \ {\bf Definition 3.1. The Sobolev-Grand Lebesgue Space} $ GW_r \psi. $ \par

\vspace{4mm}

 \ This space consists by definition from all the measurable functions having finite norm

$$
|| \ f \ ||GW_r \psi \stackrel{def}{=} \sup_{p \in (1,b)} \left\{ \frac{||f||W^{r}_p}{\psi(p)} \right\}, \
b = \const \in (1, \infty]. \eqno(3.2)
$$

\vspace{4mm}

 \ Define the following sequence of the $  \ \Psi \ -  $ functions:

 $$
 \ \theta_n(q) \stackrel{def}{=} n^{-1/q}, \ q \in (s(1),s(2)),
$$

$$
 \ s(1) = \const \in  (b, \infty), \ s(1) < s(2) = \const \le \infty,  \eqno(3.3)
$$
but the norm in this $ G\theta_n =  G\theta_n(s(1), s(2)) \  $ space, more precisely, the sequence of these norms, is defined as follows

$$
|| \ f ||G\theta_n  = || \ f ||G\theta_n(s(1), s(2)) \stackrel{def}{=}
 \sup_{q \in (s(1), s(2))} \left[ \frac{|f|_q}{\theta_n(q)} \right]. \eqno(3.4)
$$

\vspace{4mm}

 \ Denote also for brevity

$$
\Delta_n[f] = \Delta_n[f](x) :=   f(x) - J_n*f(x).
$$

\vspace{4mm}

{\bf Theorem 3.1.}

$$
|| \ \Delta_n[f] \ ||G\theta_n(s(1), s(2)) \ \le \ C_3(r) \ n^{-r} \ \frac{|| \ f \ ||GW_r \psi}{\phi(G\psi,1/n)}. \eqno(3.5)
$$

\vspace{4mm}

{\bf Proof.}  Suppose without loss of generality $  \ ||f||GW_r\psi = 1,  $ considering the value of $ r $ to be fixed; then

$$
|| \ f \ ||W^r_p \le \psi(p), 1 \le p < b.
$$

 \ We start from the well-known inequality, [1], [10], [13], [18] etc.:

$$
| \ \Delta_n[f] \ |_q \le C_3(r) \ n^{-r} \ n^{1/p - 1/q } \ ||f||W^r_p, \ 1 \le p < q, \eqno(3.6)
$$
which may be transformed as

$$
\frac{| \ \Delta_n[f] \ |_q}{n^{-1/q}} \le C_3(r) \ n^{-r} \ \frac{||f||W_p}{n^{-1/p}} \le C_3(r) \ n^{-r} \ \frac{\psi(p)}{n^{-1/p}} =
$$

$$
\frac{C_3(r) \ n^{-r} }{n^{-1/p}/\psi(p)} =  \frac{C_3(r) \ n^{-r}||f||GW\psi}{n^{-1/p}/\psi(p)}, \  1 \le p < b < q \in (s(1), s(2)). \eqno(3.7)
$$

 \ We get taking the minimum over $ \ p:  $

$$
\frac{| \ \Delta_n[f] \ |_q}{n^{-1/q}} \le C_3(r) \  n^{-r} \ \inf_{p \in (1,b)} \frac{||f||GW\psi}{n^{-1/p}/\psi(p)} =
$$

$$
C_3(r) \ n^{-r} \ \frac{||f||GW_r \psi}{\sup_p [n^{-1/p}/\psi(p)] } = C_3(r) \ n^{-r} \ \frac{||f||GW\psi}{\phi(G\psi,1/n)},\eqno(3.8)
$$
and taking further the maximum over $ \ q: $

$$
|| \ \Delta_n[f] \ ||G\theta_n(s(1), s(2)) \ \le \ C_3(r) \ n^{-r} \ \frac{||f||GW\psi}{\phi(G\psi,1/n)},
$$
Q.E.D. \\

 \ Note that the obtained statement may  be interpreted as some estimate  for Kolmogorov widths  of unit
balls in Sobolev-Grand Lebesgue Space, cf.[6], [9], [23], [13], [14], [19], [28], [33], [34]  etc.\\

\vspace{4mm}

\section{ Trigonometric approximation in Orlicz Spaces. }

\vspace{3mm}

 \ Denote by $ L(N), $   where $  N = N(u), \ u \in R $ is certain Young-Orlicz function such that

 $$
 N(u) \sim u^2, \ u \in [-1,1] \eqno(4.1)
 $$
the Orlicz function  space over source space $  X  $ with the Luxemburg norm $ || \ f \ ||L(N), \ f: X \to R.  $ \par

 \ The approximation problems by trigonometric polynomials in Orlicz spaces
were investigated by several authors, see, for example, articles [35]-[42], where was considered as a rule the case
when the generating  function $  N(u) $ satisfies the $ \Delta_2 $ condition.\par

 \ Recall that the $ \Delta_2 $ condition means the separability and reflexibility of correspondent Orlicz space. \par

 \ Let us note first of all that the so-called {\it exponential }  Orlicz space $  L(N) $ coincides up to norm equivalence
with suitable Grand Lebesgue space $  G\psi_N, $  which admit us to obtain in turn the trigonometric approximation theorems
in exponential Orlicz space. \par

\ Let $  \psi = \psi(p), \ p \in [1,b), \ b = \const, \ 1 < b \le \infty   $  (or $ p \in [1,b] ) $  be again bounded from below:
 $ \inf \psi(p) > 1 $  continuous inside the {\it  semi-open  } interval $ [1, b) $ numerical function. We can and will suppose

$$
b = \sup \{p, \ \psi(p) < \infty \},
$$
so that $ \supp \psi = [1,b)  $ or $ \supp \psi = [1,b].$ The set of all such a functions will be denoted by
$ \Psi(b), $  and we denote for brevity  $ \ \Psi:= \Psi(\infty). $ \par

 \ Suppose $ \ b = \infty $ and $  0 < ||f|| := ||f||G\psi < \infty. $ Define the function

$$
\nu(p) = \nu_{\psi}(p)  = p \ln \psi(p), \ 1 \le p < b. \eqno(4.2)
$$
 \ Recall that the Young - Fenchel, or Legendre transform $ f^*(y) $
for arbitrary function $  f: R \to R $ is defined (in the one-dimensional case) as follows

$$
f^*(y) \stackrel{def}{=} \sup_x (x y - f(x)).
$$

\ If the function $ \ f(\cdot) \ $ is continuous and convex, then

$$
f^{**}(x) = f(x),
$$
theorem of Fenchel-Moraux. \par

 \ The so - called tail function $ T_{\zeta}(y), \ \zeta: X \to R, \ y \ge 0 $ is defined by the formula

$$
T_{\zeta}(y) \stackrel{def}{=} \max \left\{ \mu(x: f(x) > y),  \mu(x: f(x)< - y)  \right\}, \ y \ge 0.
$$

 \ It is known in this case, i.e. when  $  b = \infty $  that

 $$
 T_{\zeta}(y) \le \exp \left( - \nu_{\psi}^*(\ln (y/||\zeta||) )  \right), \ y > e \cdot ||\zeta||. \eqno(4.3)
 $$
 \ Conversely, if (4.3) there holds in the following version:

$$
T_{\zeta}(y) \le \exp \left( - \nu_{\psi}^*(\ln (y/K) )  \right), \ y > e \cdot K, \ K = \const > 0,
$$
and the function $ \nu_{\zeta}(p), \ 1 \le p < \infty $ is positive, continuous, convex and such that

$$
\lim_{p \to \infty} \psi(p) = \infty,
$$
then $ \zeta \in G\psi $ and besides

$$
||\zeta||G\psi \le C(\psi) \cdot K. \eqno(4.4)
$$

\ Moreover, let us introduce the {\it exponential }  Orlicz space $ L(M) $ over the source probability space
$ (\Omega,F,{\bf P} ) $ with proper Young-Orlicz function

$$
M(u) = M_{\psi}(u):= \exp \left(  \nu_{\psi}^*(\ln |u| )  \right), \ |u| > e
$$
and as ordinary $  M(u) = M_{\psi}(u) = \exp(C \ u^2) - 1, \ |u| \le e. $ It is known \ [43] \ that the $ G\psi $ norm
of arbitrary measurable function (r.v.) $  \zeta = \zeta(x) $ is quite equivalent to the its norm in Orlicz space $ L(M): $

$$
||\zeta||G\psi \le C_1 ||\zeta||L(M) \le C_2||\zeta||G\psi, \ 1 \le C_1 \le C_2 < \infty. \eqno(4.5)
$$

 \ Evidently,  this exponential Young-Orlicz function does not satisfy the $  \ \Delta_2 \  $ condition.

\vspace{4mm}

 \ We consider now the inverse problem. Namely, let the {\it exponential} Young-Orlicz  $  M = M(u) $ be a given. The
"exponentiality"  implies by definition that the function

$$
\theta(z) = \theta_M(z) := \ln M(\exp z)
$$
is continuous and convex. Then the correspondent generating function for equivalent Grand Lebesgue Space  $  \psi_M(p) $
may be builded  by virtue of theorem Fenchel-Moraux:  $  f^{**} = f, \ f^{''} > 0, $   by the formula

$$
\psi_M(p) = \exp \left( \frac{\theta_M^*(p)}{p}  \right), \ p \ge 1. \eqno(4.6)
$$

 \ Indeed, we have

$$
 \exp \left(  \nu_{\psi}^*(\ln |u| )  \right) = M(u) = M_{\psi}(u), \ |u| > e,
$$

$$
  \nu_{\psi}^*(\ln |u| )  = \ln M_{\psi}(u), \ |u| > e,
$$
therefore

$$
  \nu_{\psi}(\ln |u| )  = \left[ \ln M_{\psi}(u) \right]^* = \theta^*_M(u),
$$
which implies (4.6). \par

\vspace{4mm}

{\bf Example 4.1.} The estimate for the r.v. $ \xi $ of a form

$$
|\xi|_p  \le C_1 \ p^{1/m} \ \ln^r p, \ p \ge 2,
$$
where $ C_1 = \const > 0, \ m = \const > 0, \ r = \const, \ $ is quite equivalent to the following tail estimate

$$
T_{\xi}(x) \le  \exp \left\{ - C_2(C_1,m,r) \ x^m \ \log^{-m r}x \right\}, \ x > e.
$$
as well as is equivalent to the belongings $ \xi(\cdot) $ to the exponential Orlicz function  with correspondent generating function
 $  N(u) = N_{m,r}(u) $  of the form

$$
N_{m,r}(u) = \exp \left( C_4(C_3,m,r) \ u^m \ \log^{-m r}u  \right), \ u > e.
$$

\vspace{3mm}

 \ It is important to note that the inequality (4.5) may be applied still when the r.v. $ \xi $ does not have the
exponential  moment, i.e. does not satisfy the famous Kramer's condition. Namely, let us consider next example. \\

\vspace{3mm}

 {\bf Example 4.2.} Define the following $ \Psi \ -  $ function.

$$
\psi_{[\beta]}(p) := \exp \left( C_3 \ p^{\beta} \right), \  p \in [1, \infty), \ \beta = \const > 0.
$$

 \ The r.v. $  \xi $ belongs to the space $ G \psi_{[\beta]} $ if and only if

$$
T_{\xi}(x) \le \exp \left( - C_4(C_3,\beta) \ [\ln (1 + x)]^{1 + 1/\beta}   \right), \ x \ge 0.
$$
as well as iff it belongs to the exponential Orlicz function  with correspondent generating function
 $  N(u) = N^{(\beta)}(u) $  of the form

$$
N^{(\beta)}(u) = \exp \left( C_5(C_3,\beta) \ln(1 + u)^{1 + 1/\beta}   \right), \ u > 1.
$$

 \ See also  [20], [24]. \par

\vspace{3mm}

 \ Let us return to the source problem. \par

\vspace{4mm}

\ {\bf Theorem 4.1.} Suppose the measurable function $ \ f: X \to R  $ belongs to certain exponential Orlicz space
$ \ L(M). \ $  We assert that this function is trigonometric approximated in this space: $  f \in TA(L(M))  $ if and only if

$$
f \in G^o \psi_M(\cdot), \eqno(4.7)
$$
or equally

$$
\lim_{\psi_M(p) \to \infty} \left\{ \frac{|f|_p}{\psi_M(p)} \right\} = 0. \eqno(4.7a)
$$

\vspace{3mm}

 \ {\bf Proof.} Since the GLS and correspondent exponential Orlicz norms ere equivalent, the problems of trigonometric
approximations in both the  considered spaces are also equal. Our proposition follows immediately from Theorem 2.1. \par

 \vspace{4mm}

 \ We recall now the following definition about comparison of Orlicz spaces. \par

 \vspace{3mm}

  \ {\bf Definition 4.1.}  Let $  L(N)  $ and $  L(K) $ be two Orlicz spaces
 with Young - Orlicz functions correspondingly $ N = N(u), \ K = K(u), \ u \ge 0.  $ We will write $ \ K << N \ $ or equally
 $ L(K) << L(N), $ iff

$$
\forall  \lambda > 0 \Rightarrow \lim_{u \to \infty} \frac{K(\lambda u)}{N(u)} = 0. \eqno(4.8)
$$

\vspace{3mm}

 \ Cf. the definition 2.1 (2.9).

 \vspace{3mm}

 \ {\bf Corollary 4.1.}  Suppose the measurable function $ \ f: X \to R  $ belongs to certain exponential Orlicz space
$ \ L(M). \ $  We assert that this function is trigonometric approximated in this space: $  f \in TA(L(M))  $ if and only if
it belongs to some Orlicz space $  \ L(K): \ f \in L(K) \  $ such that $  \ K << M. $ \par

\vspace{4mm}

\section{Concluding remarks.}

\vspace{3mm}

  \ The  multidimensional case $ \ X = [0, 2 \pi]^d, \ $  with or without weight,
 as well as  the case of trigonometric approximation on the
 whole space $  R = R^d, $ may be investigated quite analogously, as well as the problem of
algebraic approximation [1], [9], [23], [33].\par

\vspace{4mm}

\begin{center}

{\bf REFERENCES}\\

\end{center}

\vspace{4mm}

1. {\sc Ramazan Akg\"un.} {\it Approximation by polynomial in rearrangement invariant quasi Banach function spaces.}
 Banach J. Math. Anal., {\bf  6,} \ (2012), no. 2, 113-131.

\vspace{3mm}

2. {\sc Astashkin S.V.} {\it Some new Extrapolation Estimates for the Scale of
 $ L_p \ - $  Spaces.} Funct. Anal. and Its Appl., v. 37,  $ N^o $ 3, (2003), 73-77.

\vspace{3mm}

3. {\sc Belinsky E., Dai F., Ditzian Z.} {\it Multivariate approximating averages.}
  Journal of Approximation Theory. 125(1) (2003), 85-1105.

\vspace{3mm}

4. {\sc Bennet C., Sharpley R.} {\it Interpolation of operators.} Orlando, Academic Press
  Inc., (1988).

\vspace{3mm}

5. {\sc Capone C., Fiorenza A., Krbec M.} {\it On the Extrapolation Blowups in the
   $ L_p $ Scale.} Manuscripta Math., 99(4), 1999, p. 485-507.

\vspace{3mm}

6. {\sc Dai F., Ditzian Z., Tikhonov S.} {\it Sharp Jackson inequalities.}
   Journal of Approximation. Theory, (2007), 04.015

 \vspace{3mm}

7. {\sc Nina Danelia, Vakhtang Kokilashvili.} {\it Approximation by Trigonometric Polynomials in
Subspace of Weighted Grand Lebesgue Spaces.} Bulletin of the Georgian National Academy of Sciences,
vol. 7, no. 1, 2013.

\vspace{3mm}

8. {\sc Davis H.W., Murray F.J., Weber J.K.} {\it Families of $ L_p - $ spaces with
     inductive and projective topologies.} Pacific J.Math., 1970, v. 34,
    p. 619-638.

\vspace{3mm}

9. {\sc DeVore R.A., Lorentz G.G.} {\it Constructive Approximation.}  Springer, Berlin, (1993)

\vspace{3mm}

10. {\sc Borislav R. Draganov.} {\it Estimating the rate of best trigonometric
approximation in homogeneous Banach spaces by modulus of smoothness.} Proceedings of
Inst. of Mathematics and Informatics of University of Sofia Bulgarian Academy of Science, {\bf 64}, 2013, 111-132.

\vspace{3mm}

11. {\sc A.Fiorenza.} {\it Duality and reflexivity in grand Lebesgue spaces.}
       Collectanea Mathematica (electronic version), {\bf 51}, 2, (2000), 131-148.

\vspace{3mm}

12. {\sc A. Fiorenza and G.E. Karadzhov.} {\it Grand and small Lebesgue spaces and
       their analogs.} Consiglio Nationale Delle Ricerche, Instituto per le
      Applicazioni del Calcoto Mauro Picone, Sezione di Napoli, Rapporto tecnico n., {\bf 272/03,} (2005).

\vspace{3mm}

13. {\sc R. S. Ismagilov.} {\it Diameters of sets in normed linear spaces and the approximation of
functions by trigonometric polynomials.}  Russian Mathematical Surveys, 29(3), {\bf 0169,} 1974.

\vspace{3mm}

14. {\sc D. Geller, I.Z. Pesenson.} {\it Kolmogorov and linear width of balls in Sobolev spaces on compact manifold.}
arXiv:1104.0632v2 [math.FA] 27 Apr 2012

\vspace{3mm}

15. {\sc T.Iwaniec and C. Sbordone.} {\it On the integrability of the Jacobian under
      minimal hypotheses.} Arch. Rat. Mech. Anal., {\bf 119,} \ (1992), 129-143.

\vspace{3mm}

16. {\sc T.Iwaniec, P. Koskela and J. Onninen.} {\it Mapping of finite distortion:
   Monotonicity and Continuity.}  Invent. Math., {\bf 144,} (2001), 507-531.

\vspace{3mm}

17. {\sc Jawerth B., Milman M.} {\it Extrapolation Theory with Applications.}
      Mem. Amer. Math. Soc., {\bf 440,} \ (1991)

\vspace{3mm}

18. {\sc Y. Katznelson.} {\it An Introduction to Harmonic Analysis.} John Wiley  Sons, Inc., 1968.

\vspace{3mm}

19. {\sc V. N. Konovalov.} {\it Estimates of Kolmogorov type widths for classes of differentiable
periodic functions.} Mat. Zametki, 1984, Volume 35, Issue 3, 36-38.

\vspace{3mm}

20. {\sc Kozachenko Yu. V., Ostrovsky E.I.} (1985). {\it The Banach Spaces of
      random Variables of subgaussian type.}  Theory of Probab. and Math.
      Stat., (in Russian). Kiev, KSU, {\bf 32}, 43 - 57.

\vspace{3mm}

 21. {\sc Krein S.G., Petunin Yu., and Semenov E.M.} {\it Interpolation of linear  operators.} AMS, 1982.

\vspace{3mm}

22. {\sc E.Liflyand, E. Ostrovsky and L. Sirota.}
{\it Structural properties of Bilateral Grand Lebesque Spaces.}
Turk. Journal of Math., {\bf 34,} (2010), 207-219.

\vspace{3mm}

23. {\sc G.G. Lorentz, M.V. Golitschek, Yu. Makovoz.} {\it Constructive Approximation (Advanced Problems)}. Springer, Berlin, 1996.

\vspace{3mm}

24. {\sc Ostrovsky E.I.} (1999). {\it Exponential estimations for Random Fields
     and its applications}, (in Russian). Russia, OINPE.

\vspace{3mm}

25. {\sc Ostrovsky E.I.} (2002). {\it Exact exponential estimations for random
     field maximum distribution.}  Theory Probab. Appl., {\bf 45}.  v.3,  281-286.

\vspace{3mm}

26. {\sc Ostrovsky E., Sirota L.} {\it Moment Banach Spaces: Theory and Applications.}
 HIAT Journal of Science and Engineering, Holon, Israel, v. 4, Issue 1-2, (2007), 233-262.

 \vspace{3mm}

27. {\sc Ostrovsky E., Sirota L.} {\it Nikolskii-type inequalities for rearrangement invariant spaces.}
arXiv:0804.2311v1 [math.FA] 15 Apr 2008

\vspace{3mm}

28. {\sc A. Pinkus.} {\it n-widths in Approximation Theory.} Springer, New York, 1985.

\vspace{3mm}

29. {\sc G.E. Shilov.} {\it Homogeneous rings of functions.} Uspekhi Mat. Nauk, {\bf 41,} (1951),
 91 \ - \ 137 (in Russian); English translation: Amer. Math. Soc. Transl. {\bf 8,} (1954), 393-455.

\vspace{3mm}

30. {\sc Steigenwalt M.S. and While A.J.} {\it Some function spaces related to}
    $ L_p. $  Proc. London Math. Soc., 1971,  \ {\bf 22,} p. 137-163.

\vspace{3mm}

31. {\sc S.B. Stechkin.} {\it A remark on Jackson's theorem.} Tr. Mat. Inst. Steklova, {\bf 88,} (1967),
 17 \ - \ 19 (in Russian); English translation: Proc. Steklov Inst. Math. {\bf 88}, (1967), 15-17.

\vspace{3mm}

32. {\sc V. N. Temlyakov.}  {\it Approximation of periodic functions. Computational Mathematics and
 Analysis Series.} Nova Science Publishers, Inc., Commack, NY, 1993.

\vspace{3mm}

33. {\sc A.F. Timan.} {\it Theory of Approximation of Functions of a Real Variable.} Pergamon Press, 1963.

\vspace{3mm}

34. {\sc  Tikhomirov V.M. } {\it Some questions of approximation theory. } Moscow, MSU,  1976, (in Russian).

\vspace{6mm}

35. {\sc Guven, A., Israfilov, D. M.} {\it Approximation by means of Fourier trigonometric
series in weighted Orlicz spaces.} Adv. Stud. Contemp. Math. (Kyundshang), {\bf 19}, (2009), 283-295.

\vspace{3mm}

36.  {\sc Israfilov, D. M., Akg\"un, R..} {\it  Approximation in weighted Smirnov-Orlicz classes.}
J. Math. Kyoto Univ., {\bf 46,}  (2006), 755-770.

\vspace{3mm}

37. {\sc Israfilov, D. M., Guven, A.} {\it Approximation by trigonometric polynomials in
weighted Orlicz spaces.} Studia Math., {\bf 174,} (2,) (2006), 147-168.

\vspace{3mm}

38. {\sc Israfilov, D. M., Oktay, B., Akg\"un, R.} {\it Approximation in Smirnov-Orlicz classes.}
Glas. Mat. Ser. III,  {\bf 40,} (60), (2005), 87-102.

\vspace{3mm}

39. {\sc Jafarov, S. Z.} {\it Approximation by rational functions in Smirnov-Orlicz classes.}
J. Math. Anal. Appl., {\bf 379,}  (2011), 870-877.

\vspace{3mm}

40. {\sc Jafarov, S. Z.} {\it The inverse theorem of approximation of the function in Smirnov \ - \
Orlicz classes.}  Math. Inequalities, Appl. {\bf 12}, (2012), 835-844.

\vspace{3mm}

41. {\sc Jafarov, S. Z., Mamedkhanov, J. M.} {\it On approximation by trigonometric polynomials in Orlicz spaces.}
 Georgian Math. J., {\bf 19.}  (2012), 687-695.

\vspace{3mm}

42. {\sc Jafarov, S. Z.} {\it Approximation by Fejer sums of Fourier trigonometric series in
weighted Orlicz spaces.} Hacet, J. Math. Stat., {\bf 42, } (2013),  259-268.

\vspace{3mm}

43. {\sc Ostrovsky E. and Sirota L.} {\it Vector rearrangement invariant Banach spaces
of random variables with exponential decreasing tails of distributions.} \\
arXiv:1510.04182v1 [math.PR] 14 Oct 2015

\vspace{3mm}

44. {\sc Ostrovsky E.} {\it Support of Borelian measures in separable Banach spaces.} \\
arXiv:0808.3248v1 [math.FA] 24 Aug 2008

\vspace{5mm}

\end{document}